\documentclass[12pt]{article}

\usepackage{oldgerm}
\usepackage{latexsym}
\usepackage{ams}
\usepackage{amssymb}
\usepackage{amsfonts}
\usepackage{amsmath}

\def\cl{\centerline}
\def\ni{\noindent}

\def\bs{\bigskip}
\def\cen{\begin{center}}
\def\ter{\end{center}}
\def\dsp{\displaystyle}

\def\C{\mathbb{C}}

\def\go{\textgoth}

\def\Hom{\mathop{\rm Hom}\nolimits}
\def\End{\mathop{\rm End}\nolimits}

\def\dim{\mathop{\rm dim}\nolimits}
\def\deg{\mathop{\rm deg}\nolimits}

\def\Ext{\mathop{\rm Ext}\nolimits}

\def\{{\lbrace}  
\def\}{\rbrace}   
\def\({\langle}  
\def\){\rangle}
\def\[{\lbrack} \def\]{\rbrack}
\def\l{\ell}

\def\.{\bullet}

\def\rad{\hbox{rad}}

\def\ov{\overline}

\def\wt{\widetilde}

\def\ds{\mathop{\oplus}}

\def\back/{\backslash}
\def\b/{\backslash}

\def\phii#1{\f \hskip -2pt \raise -3pt \hbox{${}_#1$}}
\def\f{\varphi}

\def\iff{\Longleftrightarrow}

\def\og{\leavevmode\raise.3ex\hbox{$\scriptscriptstyle\langle\!\langle$}}
\def\fg{\leavevmode\raise.3ex\hbox{$\scriptscriptstyle\,\rangle\!\rangle$}}

\def\Lg{\boldsymbol{L}}
\def\Mg{\boldsymbol{M}}

\begin{document}

\cl{\Large \bf Introduction \`a la Conjecture d'Alexandru} \bs \bs

\ni Rappelons quelques r\'esultats de Bernstein, Gelfand, Gelfand, Delorme, 
Beilinson, Guinzburg et Soergel.  
Soient $\go{g}$ une alg\`ebre de Lie semisimple complexe, $\go{b}$ une 
sous-alg\`ebre de Borel et $\go{h}$ une sous-alg\`ebre de Cartan contenue 
dans $\go{b}$. 
Soient $\cal O$ la cat\'egorie associ\'ee \`a ces donn\'ees 
par BGG et ${\cal O}_\rho$ la sous-cat\'egorie pleine de $\cal O$ dont les 
objets ont le caract\`ere in\-fi\-ni\-t\'e\-si\-mal g\'en\'eralis\'e du module 
trivial.  
Notons $\rho$ la demi-somme des racines positives et $W$ 
le groupe de Weyl, muni de sa fonction longueur $\l$ et de son ordre de 
Bruhat.   
\`A $w \in W$ attachons le module de Verma $M_w$ de 
plus haut poids $-w\rho-\rho$~; rappelons que $M_w$ a un unique sous-module 
maximal~; notons $L_w$ le quotient correspondant.   
Soit $P_w$ un rev\`etement projectif de 
$L_w$~; posons $P:= \ds_w P_w$, $A:= (\End_{\go{g}}P)^{op}$~; notons $A$-df 
la cat\'egorie des $A$-modules de dimension finie et $E$ l'\'equivalence 
$\Hom_{\go{g}}(P,-)$ de ${\cal O}_\rho$ sur $A$-df.   
Par abus notons encore $M_w$ et $L_w$ les images de ces objets par $E$, 
et d\'esignons par $\Mg_w$ et $\Lg_w$ leurs classes respectives dans le 
groupe de Grothendieck. 
Notons $e_w \in A$ la projection sur $P_w$. 

\ni {\bf Th\'eor\`eme 1}. On a
$\dsp M_w \simeq Ae_w \Bigm/ \sum_{x \not\leq w} A \, e_x \, A \, e_w
= Ae_w \Bigm/ \sum_{x > w} A \, e_x \, A \, e_w \, .$ 

\ni {\bf Th\'eor\`eme 2}. On a $\End_A(M_w) = \C$. 

Consid\'erons les polyn\^omes de Delorme 
$a_{x,y} := SP \ \Ext_A^\.(M_x,L_y)$
o\`u $SP$ signifie $\og$s\'erie de Poincar\'e$\fg$.  

\ni {\bf Th\'eor\`eme 3}. On a $\Lg_y = \sum_x \ a_{x,y}(-1) \ \Mg_x \, .$ \bs

\vbox{\ni {\bf Th\'eor\`eme 4}. Il existe des polyn\^omes $P_{x,y}$ tels que 

\hskip 22mm (1) \ $a_{x,y} = t^{\l(y)-\l(x)} \ P_{x,y}(t^{-2}),$ 

\hskip 22mm (2) \ $P_{x,y} \not= 0 \iff x \le y \iff P_{x,y}(0) = 1,$ 

\hskip 22mm (3) \ $P_{x,x} = 1,$ 

\hskip 22mm (4) \ $\dsp \deg P_{x,y} \ < \ \frac{\l(y)-\l(x)}{2} \ 
\mbox{ si } \ x < y$.} 

\ni {\bf Th\'eor\`eme 5}. On a 
$SP \ \Ext_A^\.(L_x,L_y) \ = \ \sum_z \ a_{z,x} \ a_{z,y} \ $. 

Voici des analogues conjecturaux des ces \'enonc\'es pour les modules de 
Harish-Chandra.   

Soient $G$ un groupe de Lie semi-simple connexe \`a centre 
fini, $K$ un sous-groupe compact maximal, $\cal H$ la cat\'egorie des modules 
de Harish-Chandra associ\'ee \`a ces donn\'ees et ${\cal H}_\rho$ la 
sous-cat\'egorie pleine de $\cal H$ dont les objets ont le caract\`ere 
infinit\'esimal g\'en\'eralis\'e du module trivial. Notons $r$ le rang 
(r\'eel) de $G$ et $Z$ l'alg\`ebre $\C[[z_1,\dots,z_r]]$.   
Soit $A$ une $Z$-alg\`ebre telle que  ${\cal H}_\rho \simeq A$-df, $A$ est 
de type fini sur $Z$, $A$ est commutative modulo son radical $R$ et $A$ est 
$R$-adiquement compl\`ete.  (De telles alg\`ebres existent et sont 
isomorphes en tant que $\C$-alg\`ebres.)   
Choisissons une sous-alg\`ebre $A_0$ de $A$ relevant 
$A/R$ et notons $\{e_i \ | \ i \in I \}$ l'ensemble (fini) des idempotents 
minimaux de $A_0$.   
Soit $L_i$ le $A$-module simple associ\'e \`a 
$i \in I$, soit $\l(i)$ la dimension projective de $L_i$ et $\le$ le plus 
petit ordre sur $I$ satisfaisant $i \le j$ chaque fois que  
$$\l(j) = \l(i) + 1 \qquad \mbox{et} \qquad 
\Ext_A^1(L_j,L_i) \not= 0.$$
Utilisons librement les analogues \'evidents des notations introduites 
dans le cadre de la cat\'egorie $\cal O$. 

\ni {\bf Conjecture 1}. On a
$Ae_i \Bigm/ \sum_{j \not\leq i} A \, e_j \, A \, e_i
= Ae_i \Bigm/ \sum_{j > i} A \, e_j \, A \, e_i \, .$ 

Notons ce module $M_i$ et posons 
$$\ov{M}_i := M_i \bigm/ \rad(\End_AM_i) \, M_i \,.$$
Cet objet ne co\"{\i}ncide pas toujours avec le module de Langlands 
correspondant.

\ni {\bf Conjecture 2}. On a $\End_A(\ov{M}_i) = \C$. 

Consid\'erons les polyn\^omes de Delorme 
$a_{ij} := SP \ \Ext_A^\.(M_i,L_j)$. 

\ni {\bf Conjecture 3}. On a $\Lg_j = \sum_i \ a_{ij}(-1) \ \ov{\Mg}_i \,.$ 

\ni {\bf Conjecture 4}. Il existe des polyn\^omes $p_{ij}$ satisfaisant (1), 
\dots, (4). 

\ni {\bf Conjecture 5}. Il existe des polyn\^omes $d_k$ tels que
$$SP \ \Ext_A^\.(L_i,L_j) \ = \ \sum_k \ d_k \ a_{ki} \ a_{kj} \, .$$

Le principal inconv\'enient de cette approche des modules de 
Harish-Chandra est que, contrairement \`a ce qui se passe pour les modules de 
BGG, rien de tout cela n'est calculable~!   
Voici un rem\`ede \`a la fois partiel et conjectural \`a ce mal.   
Supposons que $G$ et $K$ ont m\^eme rang. Dans la classification 
de Langlands $L_i$ appara\^{\i}t comme l'unique quotient simple d'un module 
induit \`a partir d'un sous-groupe parabolique $P_i$ ; soit 
\hbox{${\go{p}}_i = {\go{m}}_i \ds {\go{a}}_i \ds {\go{n}}_i$} la 
d\'ecomposition de Langlands de Lie$(P_i)$ ; posons 
$$\wt d_i := \left(1-t^2 \right)^{\dim {\go a}_i} \ ;$$
soit $\wt \l(i)$ la dimension de la $K_\C$-orbite attach\'ee \`a $i$ et 
$(\wt p_{ij})$ la famille des polyn\^omes de Kazhdan-Lusztig-Vogan~; posons 
$$\wt a_{ij}(t) = t^{\wt \l(j)-\wt \l(i)} \ \wt p_{ij}(t^{-2}) \,.$$
\ni {\bf Conjecture 5'}. On a 
$SP \ \Ext_A^\.(L_i,L_j) \ = \ \sum_k \ \wt d_k \ \wt a_{ki} \ \wt a_{kj} \, $.

\bs \hrule \bs

\ni This text and others are available at 
http://www.iecn.u-nancy.fr/$\sim$gaillard

\hfill Last update : \today

\end{document}